\documentclass{amsart}
\usepackage{amsmath}%
\usepackage{amsfonts}%
\usepackage{amssymb}%
\usepackage{graphicx}
\usepackage{epstopdf}
\usepackage[margin=1in]{geometry}
\usepackage{float}
\newtheorem{theorem}{Theorem}
\theoremstyle{plain}

\newtheorem{corollary}{Corollary}

\newtheorem{definition}{Definition}

\newtheorem{lemma}{Lemma}

\numberwithin{equation}{section}

\newcommand{\N}{{\mathbb N}}
\newcommand{\Z}{{\mathbb Z}}

\def\proof{\paragraph{Proof.}}
\newcommand{\proofend}{$\Box$\bigskip}

\begin{document}
\title[Polygonal Bicycle Paths and the Darboux Transformation]{Polygonal Bicycle Paths and the Darboux Transformation}
\author{Ian Alevy}%
\curraddr{Division of Applied Mathematics, Brown University, RI 02912}%
\email{Ian\textunderscore Alevy@brown.edu}%
\author{Emmanuel Tsukerman}
\curraddr{Department of Mathematics, University of California, Berkeley, CA 94720-3840}%
\email{e.tsukerman@berkeley.edu}%
\date{}

\begin{abstract}
A Bicycle $(n,k)$-gon is an equilateral $n$-gon whose $k$ diagonals
are of equal length. In this paper we consider periodic bicycle $(n,k)$-paths,
which are a natural variation in which the polygon is replaced with a periodic polygonal path.
\end{abstract}
\maketitle

\section{Background}

Our motivation comes from three seemingly unrelated problems. The first is the
problem of \textit{Floating Bodies of Equilibrium in Two Dimensions}.
From 1935 to 1941, mathematicians at the University of Lviv, among
them Stefan Banach and Mark Kac, collected mathematical problems in
a book, which became known as the Scottish Book, since they often
met in the Scottish Coffee House. Stanislaw Ulam posed problem 19
of this book: ``Is a sphere the only solid of uniform
density which will float in water in any position?'' The answer in the two-dimensional case, as it turns out, depends on the density of the solid.

\begin{figure}[H]
\centering
\includegraphics[width=4in]{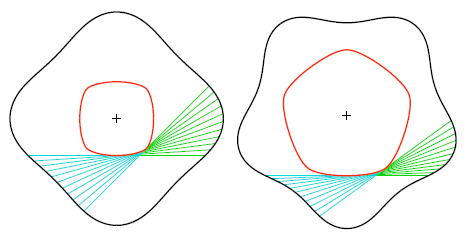}
\caption{\label{fig:Ambitrack}Ambiguous bicycle tracks. The rear-wheel track is the inner curve and the front-wheel track is the outer curve. One cannot tell which way the bicycle went because a bicycle could have followed either one of two  trajectories \cite{Wegner1}. }
\end{figure}

The second problem, known as the \textit{Tire Track Problem}, originated
in the story, ``The Adventure of the Priory School''
by Arthur Conan Doyle, where Sherlock Holmes and Dr. Watson discuss
in view of the two tire tracks of a bicycle which way the bicycle
went. The problem is: ``Is it possible that tire
tracks other than circles or straight lines are created by bicyclists
going in both directions?'' As shown in Figure \ref{fig:Ambitrack}, the answer to this subtle question is affirmative.

The third problem is the problem of describing the trajectories of
\textit{Electrons in a Parabolic Magnetic Field}. All three problems
turn out to be equivalent \cite{Wegner1}. 

Often in mathematics it is fruitful to discretize
a problem. As such, S. Tabachnikov proposed a ``discrete bicycle
curve'' (also known as a ``bicycle polygon'') \cite{Tabachnikov1}, which is a polygon satisfying discrete analogs
of the properties of a bicycle track. The main requirement turns out
to be that, in the language of discrete differential geometry, the
polygon is ``self-Darboux''. That is, the discrete differential
geometric notion of a Discrete Darboux Transformation \cite{DDG1, DDG2}, which
relates one polygon to another, relates a discrete bicycle curve to
itself. 

The topic of bicycle curves and polygons belongs to a number of active areas of research. On the one hand, it is part of rigidity theory. As an illustration, R. Connelly and B.  Csik\'{o}s consider the problem of classifying first-order flexible regular bicycle polygons \cite{ConnellyCsikos}. Other work on the rigidity theory of bicycle curves and polygons can be found in \cite{Csikos, Cyr, Tabachnikov1}.

The topic is also part of the subject of Discrete Integrable Systems. This point of view is taken in \cite{TabTsu}, where the authors find integrals of motion (i.e., quantities which are conserved) of bicycle curves and polygons under the Darboux Transformation and Recutting of polygons \cite{Adler1, Adler2}.

In this paper, in analogy with bicycle polygons, we introduce a new concept called a periodic discrete bicycle path and study both its rigidity and integrals.

\section{Bicycle $(n,k)$-paths}
 A Bicycle $(n,k)$-gon is an equilateral $n$-gon whose $k$ diagonals
are of equal length \cite{Tabachnikov1}. We consider the following analog.

\begin{definition} \label{defBikePath} 
Define $P=\{V_{i}\in\mathbb{R}^{2}:i\in\mathbb{Z}\}$
(for brevity, $V_{0}V_{1}\cdot\cdot\cdot V_{n-1}$) to be a \textit{discrete
periodic bicycle $(n,k)$-path }(or discrete $(n,k)$-path) if the
following conditions hold:
\begin{enumerate}  
\item[(i).]$V_{n+i}=V_{i}+e_{1}\mbox{ }\forall i$, $V_{0}=(0,0)$. (Periodicity Condition)
\item[(ii).] $|V_{i}V_{i+1}|=|V_{j}V_{j+1}|\forall i,j$. (Equilateralness)
\item[(iii).] $|V_{i}V_{i+k}|=|V_{j}V_{j+k}|\forall i,j$. (Equality of $k$-diagonals)
\end{enumerate}
\end{definition}

\begin{figure}[H]
\centering
\includegraphics[width=4in]{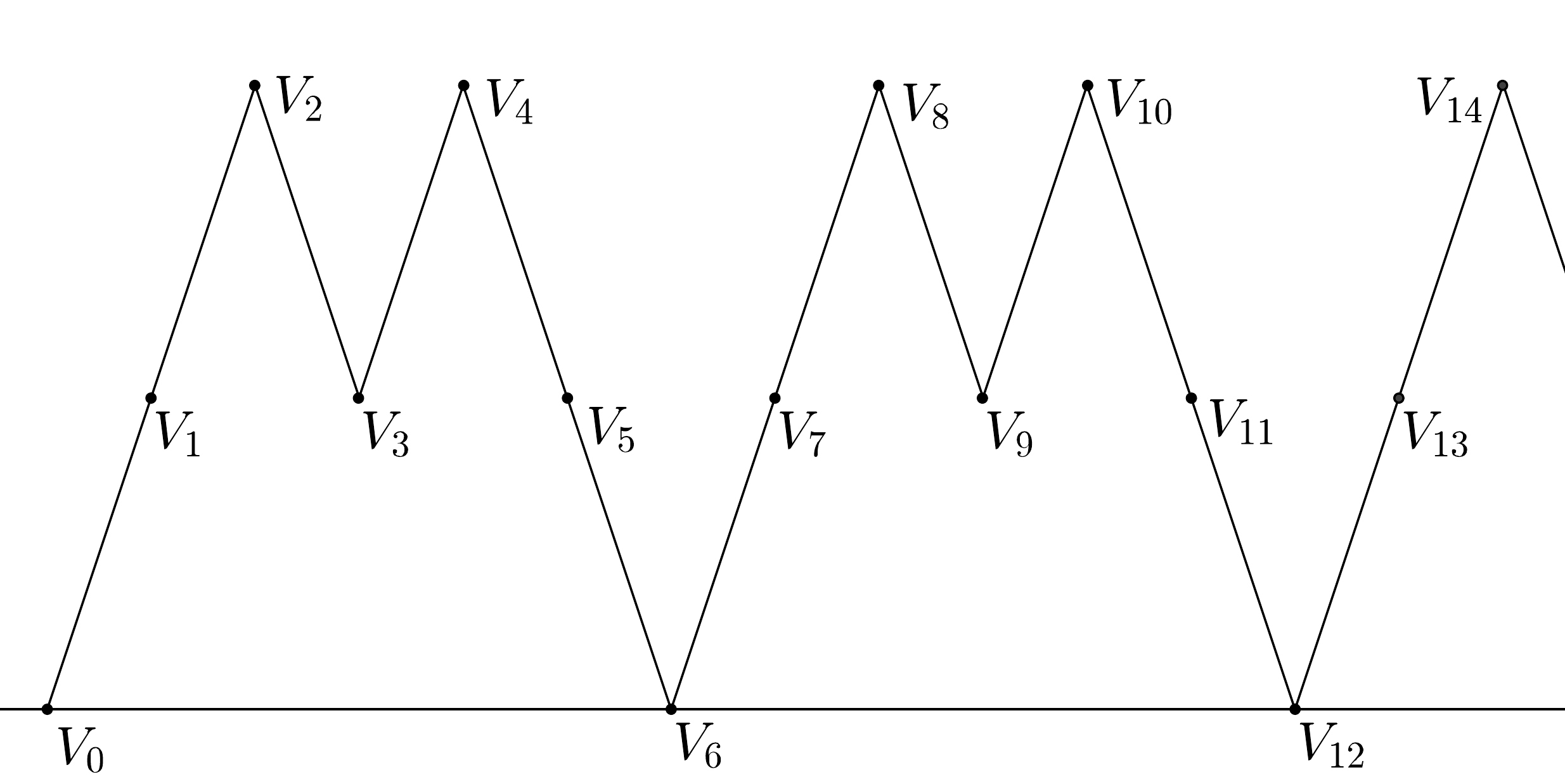}
\caption{\label{sixFive} An example of a discrete $(6,5)$-path.}
\end{figure}

Definition \ref{defBikePath} is meant to model the motion of a bicycle whose trajectory is spatially periodic. If the back wheel of the bike is at $V_j$ then the front wheel is at $V_{j+k}$, and the condition that $|V_j V_{j+k}|$ is independent of $j$ represents the rigidity of the frame of the bike. The condition that $|V_j V_{j+1}|$ is independent of $j$ prescribes a constant speed for the motion of the bike.

Some natural questions regarding periodic $(n,k)$-paths are for which pairs $(n,k)$ they exist, how many there are and whether they are rigid or flexible. We consider these questions in section \ref{sec:rig}. A simple example of a bicycle $(n,k)$-path, analogous to the regular $(n,k)$-polygon, is $V_{i}=(\frac{i}{n},0)$, i.e. when all vertices lie at equal intervals on the line. We call this the \textit{regular path}. 
Since bicycle $(n,k)$-paths are discretized bicycle paths, it is also interesting to see if there are any integrals of motion. We show that this is indeed the case in section \ref{sec:integrals}, by showing that area is an integral of motion.

\section{Rigidity \label{sec:rig}}

The following two Lemmas will be helpful when analyzing the rigidity of discrete bicycle paths.
\begin{lemma}
\label{lem:cyclic absolute values}Let $n\in\mathbb{N}$, $\chi_{i}\in\{-1,1\}$ for every $i\in {\Z} / n \Z$
and let 
\[
S=\{(x_{0},x_{1},...,x_{n-1})\in\mathbb{R}^{n}:(x_{i+1}-x_{i})^{2}=(x_{j+1}-x_{j})^{2}\mbox{ }\forall i,j\in\Z / n \Z\}.
\]
Then 
\[
S=\{(x_{0},x_{1},...,x_{n-1}):x_{i+1}=x_{i}+\chi_{i}r\mbox{ for }i\in\Z / n \Z,\sum_{i=0}^{n-1}r\chi_{i}=0\mbox{ and }r\geq0\}.
\]
 In particular, if $n$ is odd, then $S=\{(x_{0},x_{1},...,x_{n-1}):x_{i}=x_{j}\mbox{ }\forall i,j\in\Z / n \Z\}$. \end{lemma}
\proof
First note that the candidate set is well-defined since $x_{j+n}=x_{j}+\sum_{i=j}^{j+n-1}r\chi_{i}=x_{j}+\sum_{i=0}^{n-1}r\chi_{i}=x_{j}$.
Let $(x_0,x_1,\ldots,x_{n-1}) \in  S$. Recall that
\[
\mbox{sgn}(x)=\begin{cases}
1 & \mbox{if }x\geq0\\
-1 & \mbox{if }x<0
\end{cases}
\]
 and that $\mbox{sgn}(x)|x|=x$. Set $r:=|x_{i+1}-x_{i}|$ and $\chi_{i}=\mbox{sgn}(x_{i+1}-x_{i})$. Then
\[
x_{i+1}=x_{i}+\chi_{i}r
\]
and
\[
\sum_{i=0}^{n-1}r\mbox{sgn}(x_{i+1}-x_{i})=0.
\]
It follows that any $n$-tuple in $S$ satisfies the conditions $x_{i+1}=x_{i}+\chi_{i}r$, $\sum_{i=0}^{n-1}r\chi_{i}=0$
and $r\geq0$. The opposite inclusion is clear.
\proofend

\begin{lemma}
\label{lem:difference equation}Let $x_{i}\in\mathbb{R}$ for every
$i\in\mathbb{Z}$ with $x_{0}=0$ and let $k$ and $n$ be coprime
integers. Assume that $x_{i+k}-x_{i}=x_{i}-x_{i-k}$ for each $i$
and that $x_{i+n}=1+x_{i}$. Then $x_{i}=\frac{i}{n}$ for each $i$.\end{lemma}
\proof 
Define $z_i$ via $x_i=z_i+\frac{i}{n}$. The hypothesis $x_{i+n}=1+x_{i}$ implies that
\[
z_{i+n}=z_i.
\]
Since $k$ and $n$ are coprime, the difference
\[
\Delta z:=z_{i+k}-z_i
\]
is independent of $i$.
If $|k|$ is the order of $k$ in $\Z / n \Z$, then
\[
z_{i+n}-z_i=z_{i+|k|k}-z_i=|k| \Delta z = 0,
\]
so that $z_i=0$ for every $i$.
\proofend

The following Theorem gives a classification of a family of periodic $(n,k)$-paths.
\begin{theorem} \label{Thm:discretePath1} The discrete $(n,dn-1)$-paths $V_i=(x_i,y_i), i \in \N$ with $d\neq0$ are exactly those paths which satisfy 
\[
x_{j}=\frac{j}{n}
\]
 and \textup{
\[
y_{j+1}=y_{j}+\chi_{j}r\mbox{ for }j\in\Z / n \Z\mbox{ with }\sum_{i=0}^{n-1}r\chi_{i}=0\mbox{ and }r\geq0
\]
}
for each $j$.
In particular, if $n$ is odd then a discrete $(n,dn-1)$-path must
be regular.\end{theorem}
\proof
For every $i$,
\[
|V_{i}V_{i+1}|=|V_{i+dn-1}V_{i+dn}|
\]
and 
\[
|V_{i}V_{i+dn-1}|=|V_{i+1}V_{i+dn}|
\]
Therefore
\[
(x_{i+1}-x_{i})^{2}+(y_{i+1}-y_{i})^{2}=(x_{i}-x_{i-1})^{2}+(y_{i}-y_{i-1})^{2}
\]
and
\[
(d+x_{i-1}-x_{i})^{2}+(y_{i-1}-y_{i})^{2}=(d+x_{i}-x_{i+1})^{2}+(y_{i}-y_{i+1})^{2}.
\]
It follows that 
\[
d(x_{i+1}-x_{i})=d(x_{i}-x_{i-1})
\]
Since $d\neq0$, 
\[
x_{i+1}-x_{i}=x_{i}-x_{i-1}.
\]
By Lemma \ref{lem:difference equation}, $x_{j}=\frac{j}{n}$ for
each $j$. Now equation $|V_{i}V_{i+1}|=|V_{j}V_{j+1}|$ for all $i,j$
implies that 
\[
(x_{i+1}-x_{i})^{2}+(y_{i+1}-y_{i})^{2}=(x_{j+1}-x_{j})^{2}+(y_{j+1}-y_{j})^{2}\mbox{ }\forall i,j
\]
so that
\[
(y_{i+1}-y_{i})^{2}=(y_{j+1}-y_{j})^{2}
\]
By Lemma \ref{lem:cyclic absolute values}, we are done.
\proofend

\begin{theorem}
\label{thm:path2} The discrete $(n,dn+1)$-paths $V_i=(x_i,y_i), i \in \N$ with $d\neq0$ are exactly those paths which satisfy 
\[
x_{j}=\frac{j}{n}
\]
 and \textup{
\[
y_{j+1}=y_{j}+\chi_{j}r\mbox{ for }j\in\Z / n \Z\mbox{ with }\sum_{i=0}^{n-1}r\chi_{i}=0\mbox{ and }r\geq0
\]
}
for each $j$.
In particular, if $n$ is odd then a discrete $(n,dn+1)$-path must
be regular.
\end{theorem}
\proof
Set $C_{1}=|V_{i}V_{i+dn+1}|^{2}$ and $C_{2}=|V_{i}V_{i+1}|^{2}$.
Then
\[
(d+x_{i+1}-x_{i})^{2}+(y_{i+1}-y_{i})^{2}=C_{1}
\]
and 
\[
(x_{i+1}-x_{i})^{2}+(y_{i+1}-y_{i})^{2}=C_{2}.
\]
Substituting, we get
\[
d^{2}+2d(x_{i+1}-x_{i})+C_{2}=C_{1}
\]
so that $x_{i+1}-x_{i}$ is constant. By Lemma \ref{lem:difference equation}, $x_{i}=\frac{i}{n}$. It follows that $(y_{i+1}-y_{i})^{2}$ is constant, so by Lemma \ref{lem:difference equation} we are done.
\proofend

\begin{corollary}
Any $(n,dn+1)$-path is an $(n,dn-1)$-path and vice verse.
\end{corollary}

For an example, see Figure \ref{sixFive}.

\section{Darboux Transformation and Integrals \label{sec:integrals}}

It is important to make a distinction between infinitesimal ``trapezoidal''
movement and infinitesimal ``parallelogram'' movement of the bicycle. Informally, infinitesimal ``parallelogram'' movement can be seen as a trivial gliding of the bike, and thus we insist on infinitesimal ``trapezoidal'' movement.
\begin{definition}
(Trapezoidal Condition) We will say that a discrete $(n,k)$-path
satisfies the trapezoidal condition if $V_{i}V_{i+k+1}\parallel V_{i+1}V_{i+k}$
for each $i\in \Z$.
\end{definition}

 As an illustration of these concepts, consider Figure \ref{sixFive}: $V_0 V_1 V_5 V_6$ is a trapezoidal motion, while $V_1 V_2 V_6 V_7$ is a parallelogram motion. Consequently, the bicycle path in the figure does not satisfy the Trapezoidal Condition.

Assuming the Trapezoidal Condition, we may view bicycle paths
in terms of an important construction in Discrete Differential Geometry
called the Darboux Transformation \cite{DDG1, DDG2}.
\begin{definition}
(Darboux Transform) We say that two polygons $P$ and $Q$ are in Darboux Correspondence (with parameter $\ell$) if each quadrilateral
$P_{i}Q_{i}P_{i+1}Q_{i+1}$ is an isosceles trapezoid with side length
$\ell$:
\[
|P_{i}Q_{i}|=|P_{i+1}Q_{i+1}|=\ell \mbox{ and }P_{i}Q_{i+1}\parallel Q_{i}P_{i+1}.
\]
\end{definition}

\begin{figure}[H]
\centering
\includegraphics[width=3in]{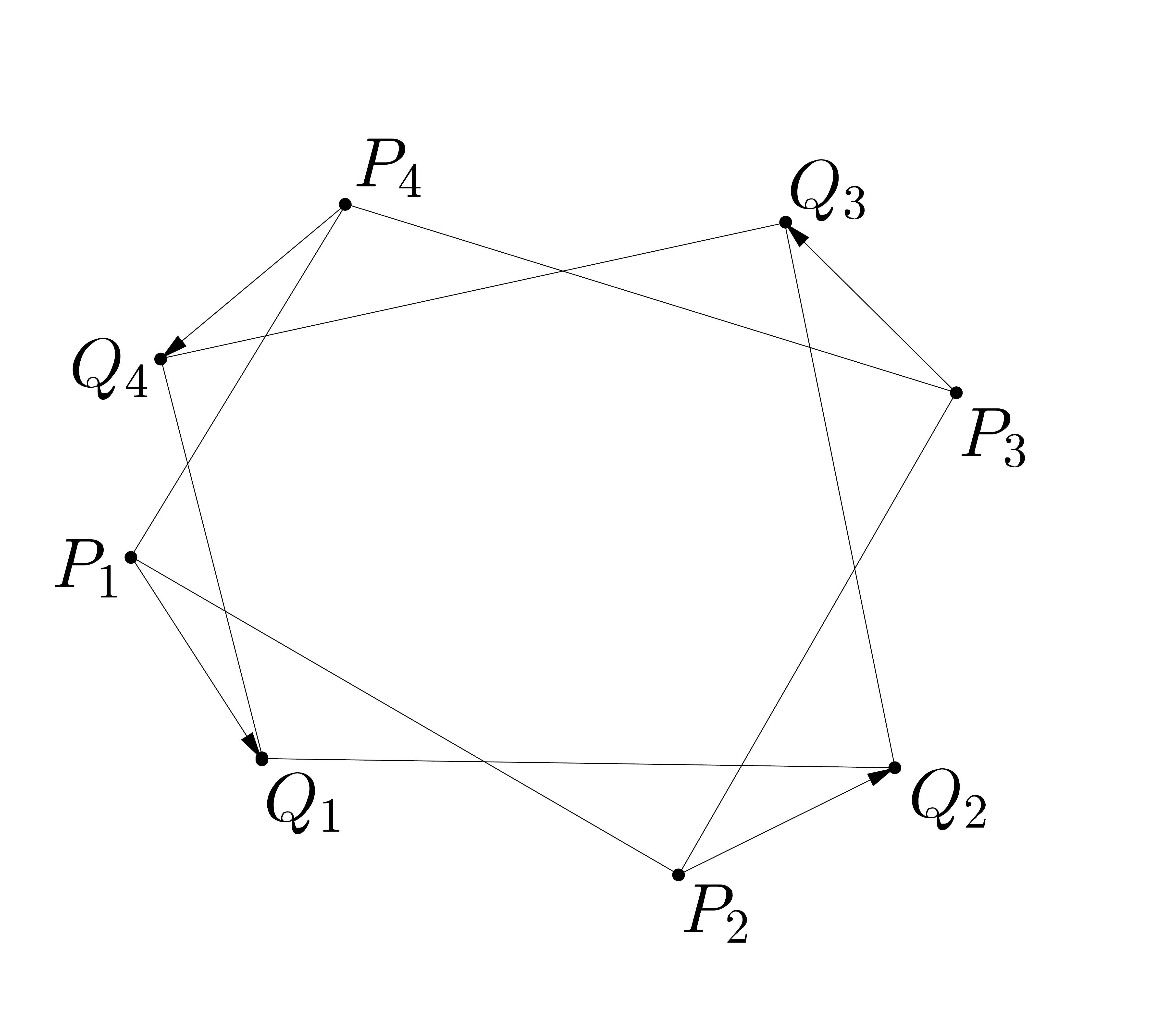}
\caption{\label{DCorrespondence}Two polygons in Darboux Correspondence. }
\end{figure}

We denote the map taking vertex $P_i$ to $Q_i$ by
$\mathfrak{D}$. We will also refer to the map of polygons $\mathfrak{D}(P)=Q$ by the same letter, since no confusion ought to occur.

Consider a polygonal line $P$ with vertices $V_{0},V_{1},...,V_{n-1}$.
Let $v_{0}$ be a vector with its origin at $V_{0}$. Having a vector
$v_{i}$ at vertex $V_{i}$, we obtain a vertex $v_{i+1}$ of the
same length at $V_{i+1}$ via the trapezoidal condition. For example, in Figure \ref{DCorrespondence}, $v_1=P_1 Q_1$ and $v_2=P_2 Q_2$. For a fixed length
of $v_{0}$, we may view the map taking $v_{0}$ to $v_{j}$ as a
self-map of the circle of radius $|v_{0}|=|v_{j}|$ by identifying
the circle at $V_{0}$ with circle at $V_{j}$ via parallel translation. 
\begin{definition}
(Monodromy map of the Darboux Transformation) The monodromy map is the map acting on the identified
circles at $V_{0}$ and $V_{n}$ which takes $v_{0}$ to $v_{n}$. 
\end{definition}
It is known  that the Monodromy map is a fractional linear transformation
on a circle of fixed radius after we identify the circle with the real projective line via stereographic projection \cite{TabTsu}. We will assume throughout, unless otherwise
stated, that the monodromy map is acting on a fixed point. In other
words, we will assume that the Darboux transform has been chosen so
that the Darboux vector at $V_{0}$ is the same as that at $V_{n}$,
where $n$ is the period. This is analogous to applying the Darboux
transform to a closed polygon and requiring that its image is closed
also. \\

We mention in passing that in the case of closed polygons, Darboux Correspondence implies that the monodromies of the two polygons are conjugated to each other. The invariants of the conjugacy class of the monodromy, viewed as functions of the length parameter are consequently integrals of the Darboux Correspondence \cite{TabTsu}.

\begin{figure}[H]
\includegraphics[scale=0.36]{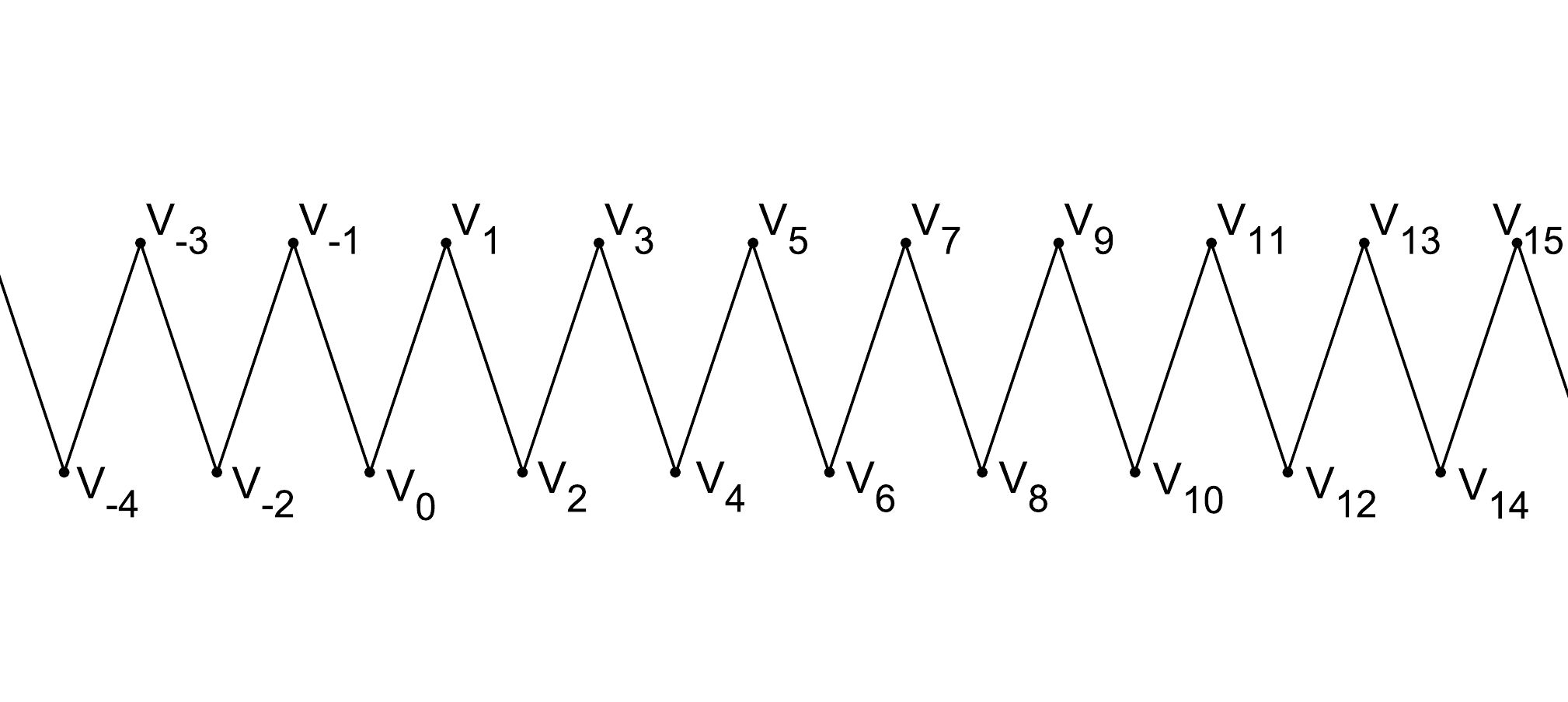} \quad \includegraphics[scale=0.36]{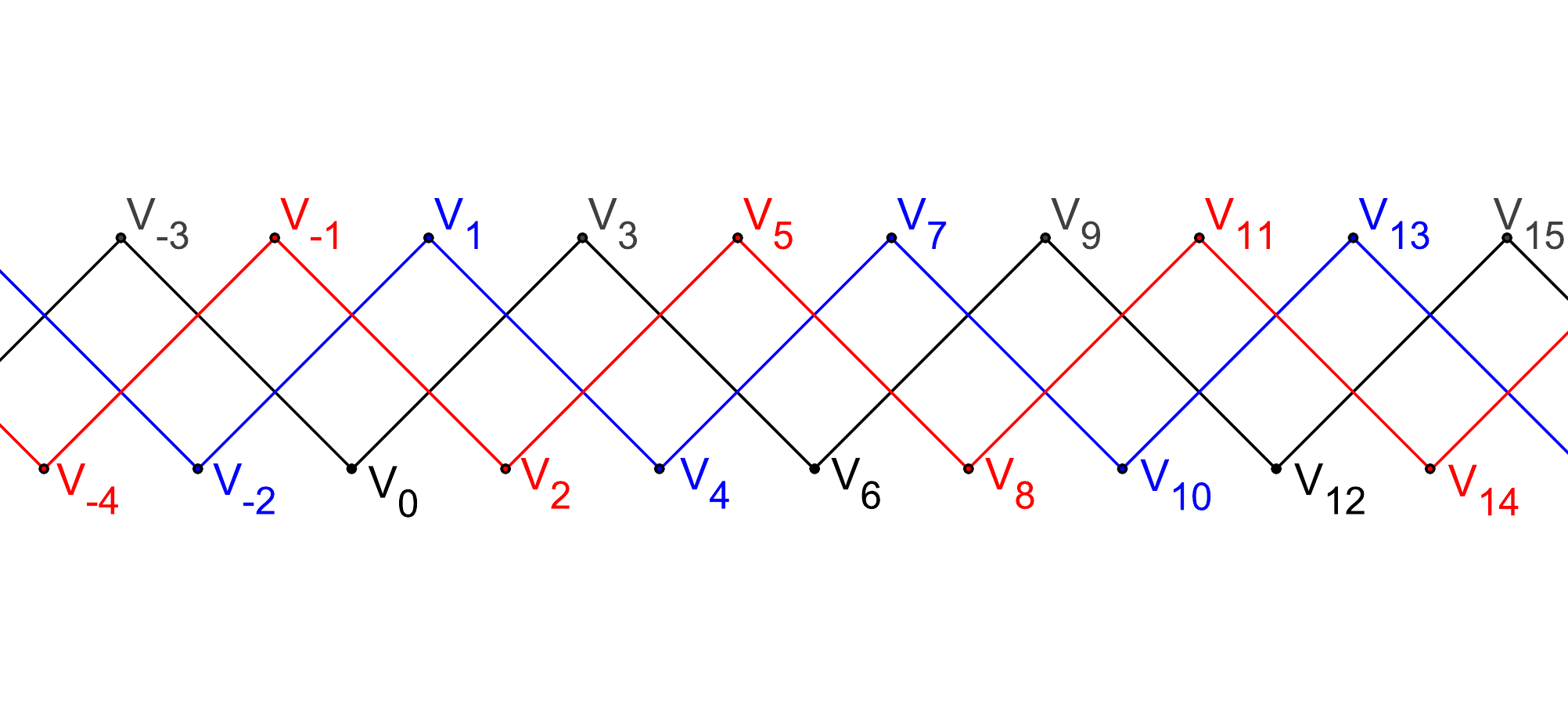}

\includegraphics[scale=0.36]{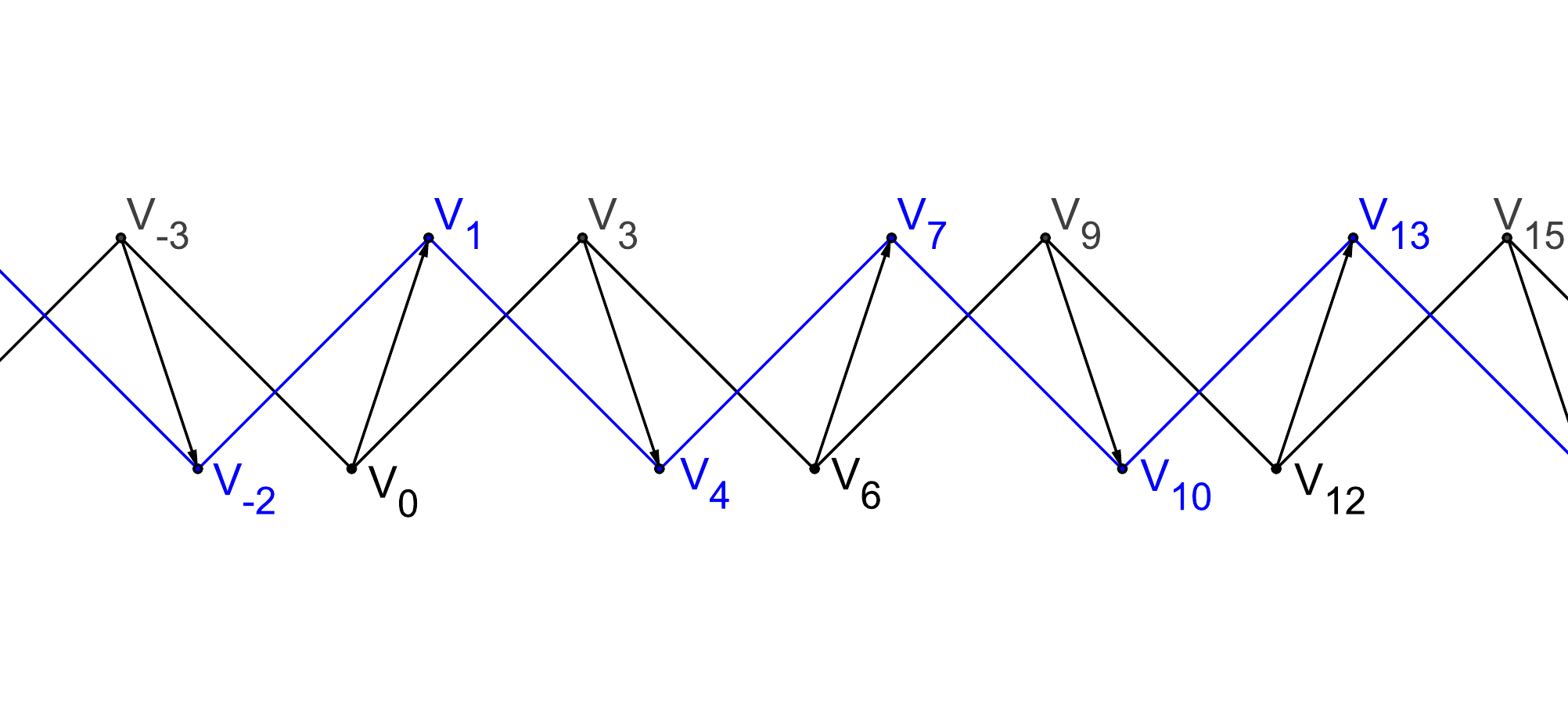}

\caption{ \label{connection} Viewing a discrete $(n,k)$-path satisfying the Trapezoidal Condition (top left) in terms of the Darboux Transformation. The path is decomposed into equilateral linkages (top right). Any two consecutive linkages are in Darboux Correspondence (bottom).}

\end{figure}

\textbf{Connection between Darboux Transformation and Discrete $(n,k)$-paths}. A discrete $(n,k)$-path satisfying the trapezoidal condition may
be interpreted in terms of the Darboux Transform. Indeed, given such a path, we consider the periodic equilateral linkages $L_i=\ldots V_{0+i} V_{k+i} V_{2k+i} \ldots$ for $i=0,1,\ldots,k-1$. The trapezoidal condition implies that there is a Darboux Correspondence $\mathfrak{D}(L_i)=L_{i+1}$ of the same parameter (since the $(n,k)$-path is equilateral) for consecutive linkages (see Figure \ref{connection}).
Conversely, given a series of such equilateral linkages in consecutive Darboux Correspondence, it is easy to see that we may form a bicycle $(n,k)$-path satisfying the trapezoidal condition.  \\

\begin{figure}[hbtp]
\centering
\includegraphics[width=4in]{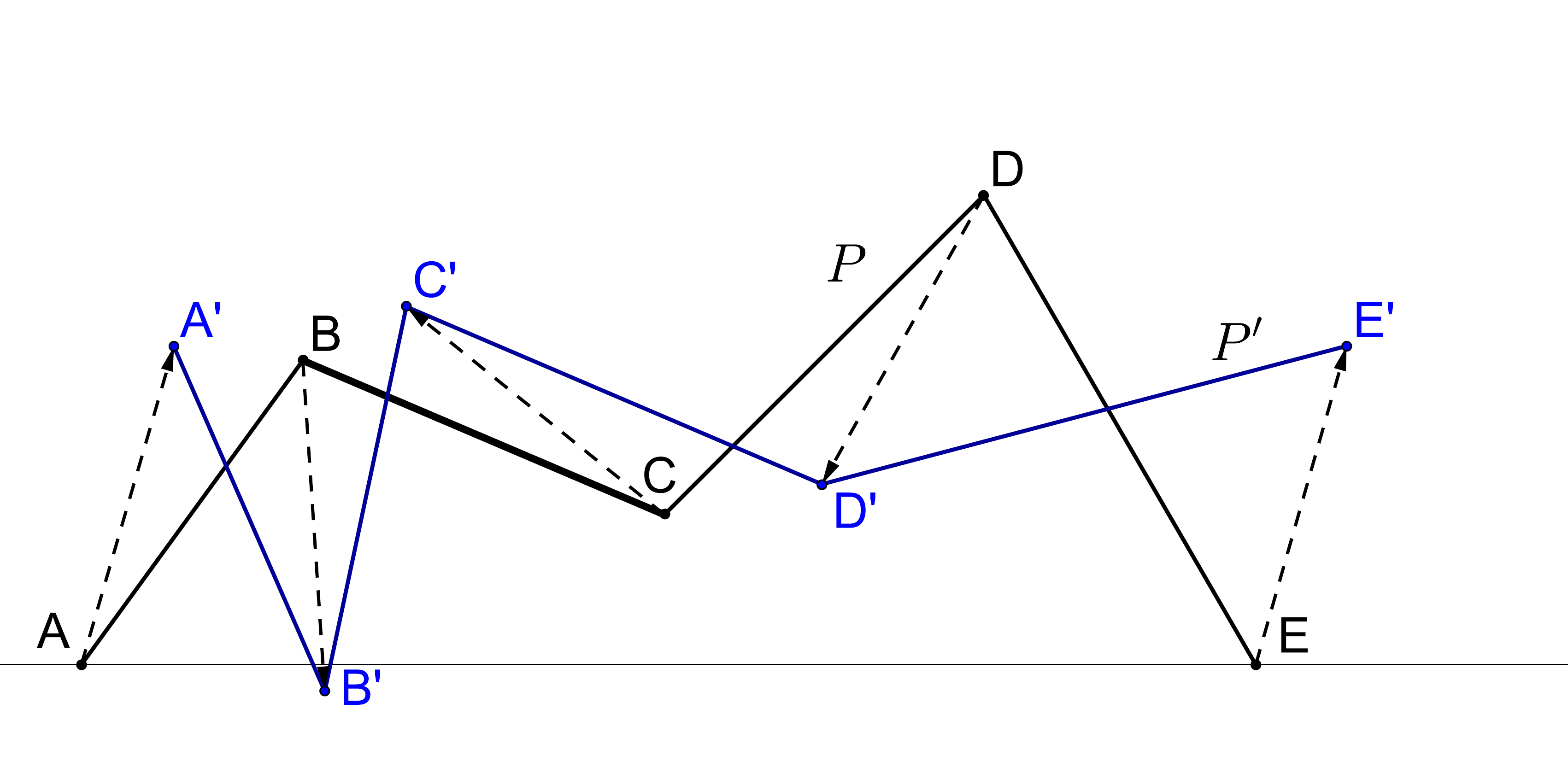}
\includegraphics[width=4in]{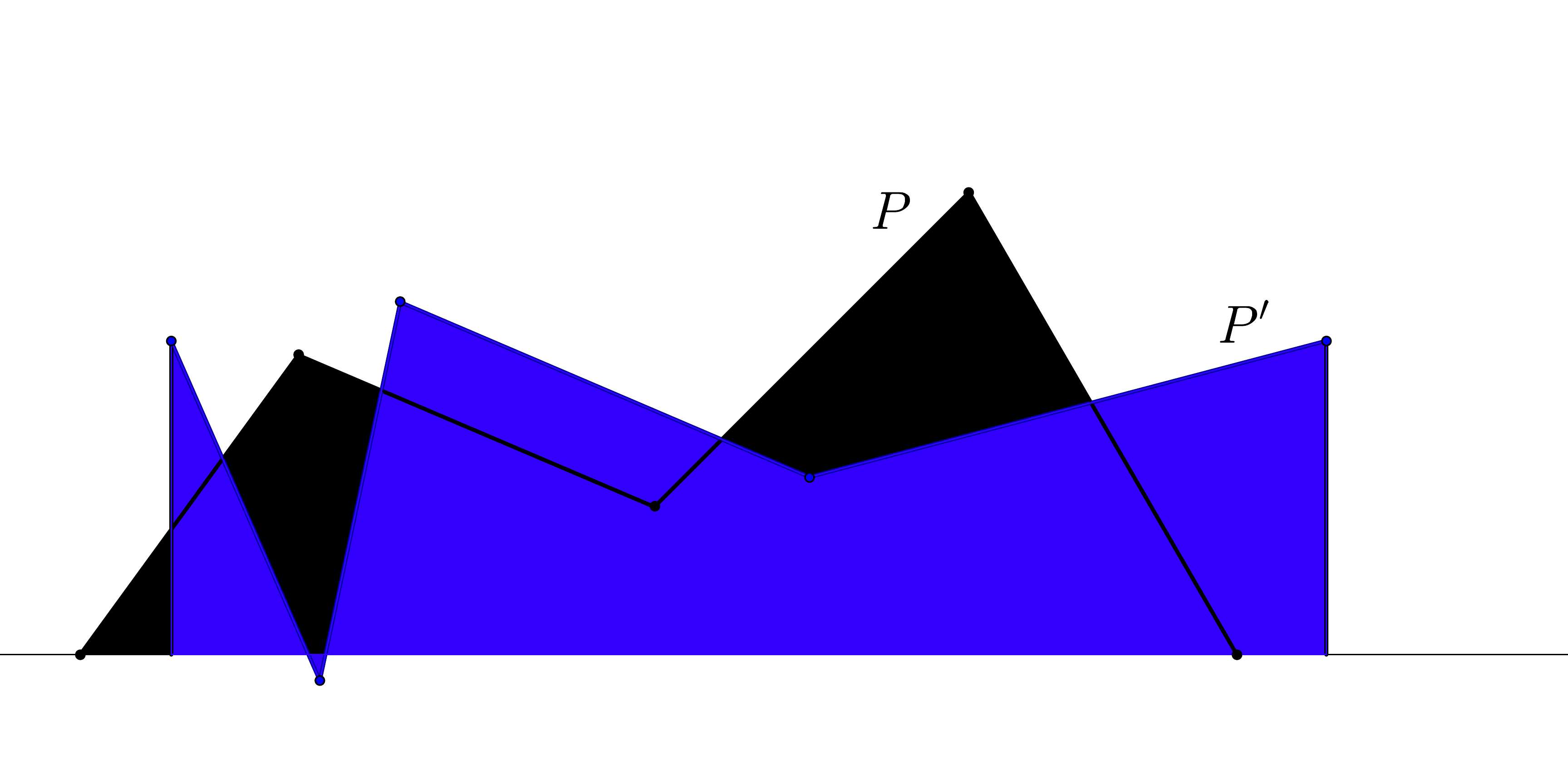}
\caption{\label{areaFig}Two periodic paths $P$ and $P'$ in Darboux Correspondence. By Theorem \ref{thm:area}, the two paths have equal areas under the curve.}
\end{figure}

The Darboux Transformation also preserves the area of periodic linkages.
More precisely, let $y=-c$ for $c>0$ sufficiently large so that
the linkage $L$ and its Darboux Transformation $D$ lie completely
above $y=-c$. We define an area function as follows. Let
$\check{V_{i}}=(x(V_{i}),-c)$. We define the area of $P$ to be the
signed area of the polygon $\check{V_{0}}V_{0}V_{1}\cdot\cdot\cdot V_{n}\check{V_{n}}$
and denote it by $|P|$. We show that this area is preserved under Darboux Transformation (see Figure \ref{areaFig}). In particular, it will follow that 
the area of $V_{0}V_{k}\cdot\cdot\cdot V_{n-k}$ is equal to the area
of $V_{m}V_{k+m}\cdot\cdot\cdot V_{n-k+m}$ for every $m\in\mathbb{Z}$.

\begin{theorem} \label{thm:area}
The Darboux transformation is area preserving on periodic polygonal paths.
\end{theorem}

\proof
Let $P$ and $D$ be two periodic polygonal paths in Darboux correspondence. We show that the difference of the areas of $P$ and $D$ is zero.
We have 
\[
|P|=\sum_{i=0}^{n-1}|\check{V}_{i}V_{i}V_{i+1}\check{V}_{i+1}|
\]

and similarly for $D$. Therefore
\[
|P|-|D|=\sum_{i=0}^{n-1}|\check{V}_{i}V_{i}V_{i+1}\check{V}_{i+1}|-|\check{V}_{i}'V_{i}'V_{i+1}'\check{V}_{i+1}'|.
\]

From the isosceles trapezoids, 
\begin{equation}
|V_{i}V_{i+1}V_{i+1}'|=|V_{i}'V_{i+1}'V_{i}|.\label{eq:Trapezoid Triangle Area Equality}
\end{equation}

Also,

\[
|\check{V}_{i}V_{i}V_{i+1}\check{V}_{i+1}|=|\check{V}_{i}V_{i}V_{i+1}'\check{V_{i+1}'}|+|\check{V}_{i+1}'V_{i+1}'V_{i+1}\check{V_{i+1}}|+|V_{i}V_{i+1}V_{i+1}'|.
\]

Similarly, 
\[
|\check{V}_{i}'V_{i}'V_{i+1}'\check{V}_{i+1}'|=|\check{V}_{i}'V_{i}'V_{i}\check{V_{i}}|+|\check{V}_{i}V_{i}V_{i+1}'\check{V_{i+1}}'|+|V_{i}'V_{i+1}'V_{i}|
\]

Using equation \ref{eq:Trapezoid Triangle Area Equality},
\[
|\check{V}_{i}V_{i}V_{i+1}\check{V}_{i+1}|-|\check{V}_{i}'V_{i}'V_{i+1}'\check{V}_{i+1}'|=|\check{V}_{i+1}'V_{i+1}'V_{i+1}\check{V_{i+1}}|-|\check{V}_{i}'V_{i}'V_{i}\check{V_{i}}|.
\]

It follows that 
\[
|P|-|D|=\sum_{i=0}^{n-1}|\check{V}_{i+1}'V_{i+1}'V_{i+1}\check{V_{i+1}}|-|\check{V}_{i}'V_{i}'V_{i}\check{V_{i}}|,
\]

which telescopes to 
\[
|P|-|D|=|\check{V}_{n}'V_{n}'V_{n}\check{V_{n}}|-|\check{V}_{0}'V_{0}'V_{0}\check{V_{0}}|.
\]

Since $\vec{V_{n}V_{n}'}=\vec{V_{0}V_{0}'}+e_{1}$ and $V_{n}=V_{0}+e_{1}$,
it follows that $|\check{V}_{n}'V_{n}'V_{n}\check{V_{n}}|=|\check{V}_{0}'V_{0}'V_{0}\check{V_{0}}|$,
so that $|P|=|D|$.
\proofend

\section{Questions}

We end our discussion with some research topics and questions of interest concerning bicycle $(n,k)$-paths.

\begin{enumerate}
\item Construct interesting families of bicycle $(n,k)$-paths. For example, ones for which the condition $x_j=\frac{j}{n}$ does not hold.
\item What is the $m$th order infinitesimal rigidity theory of bicycle $(n,k)$-paths like?
\item For closed bicycle polygons, there are many integrals of motion \cite{TabTsu}. For example, a geometric center called the Circumcenter of Mass \cite{TabTsu2} is invariant under Darboux Transformation for closed polygons. Are there other integrals of motion for bicycle $(n,k)$-paths?
\end{enumerate}

\bigskip
{\bf Acknowledgments}. It is a pleasure to acknowledge the valuable help and discussion with S. Tabachnikov. This project originated during the program Summer@ICERM 2012; we are grateful to ICERM for the support.

\end{document}